\documentclass[twoside,a4paper]{article}

\usepackage{amsmath,amsthm,amssymb,latexsym}
\usepackage[v2,tips]{xy}

\newcommand{\ip}[2]{{\langle {#1} , {#2} \rangle}}

\newcommand{\mc}[1]{\mathcal{#1}}

\newcommand{\inten}{{\check{\otimes}}}
\newcommand{\proten}{{\widehat{\otimes}}}

\newcommand{\exthaa}{\otimes^{eh}}
\newcommand{\norhaa}{\otimes^{\sigma h}}

\theoremstyle{plain}%
\newtheorem{proposition}{Proposition}[section]%
\newtheorem{theorem}[proposition]{Theorem}%
\newtheorem{corollary}[proposition]{Corollary}%
\newtheorem{lemma}[proposition]{Lemma}%
\theoremstyle{definition}%
\newtheorem{definition}[proposition]{Definition}%
\theoremstyle{remark}%

\begin{document}

\large
\title{\textsc{Conditions implying the uniqueness of the weak$^*$-topology on certain group algebras}}
\author{Matthew Daws\\\normalsize\texttt{matt.daws@cantab.net}\and Hung Le Pham\thanks{Supported by a Killam postdoctoral fellowship}\\\normalsize\texttt{hlpham@math.ualberta.ca }\and Stuart White\\\normalsize\texttt{s.white@maths.gla.ac.uk}}
\maketitle

\begin{abstract}
We investigate possible preduals of the measure algebra $M(G)$ of a locally compact group and the Fourier algebra $A(G)$ of a separable compact group.  Both of these algebras are canonically dual spaces and the canonical preduals make the multiplication separately weak$^*$-continuous so that these algebras are dual Banach algebras.  In this paper we find additional conditions under which the preduals $C_0(G)$ of $M(G)$ and $C^*(G)$ of $A(G)$ are uniquely determined. In both cases we consider a natural coassociative multiplication and show that the canonical predual gives rise to the unique weak$^*$-topology making both the multiplication separately weak$^*$-continuous and the coassociative multiplication weak$^*$-continuous. In particular, dual cohomological properties of these algebras are well defined with this additional structure.
\end{abstract}

\section{Introduction}
A classical theorem of Sakai characterises von Neumann algebras as those $C^*$-algebras which are isometrically isomorphic to the dual space of some Banach space.  Furthermore, the weak$^*$-topology induced by this duality is the unique weak$^*$-topology making the multiplication separately continuous and the adjoint continuous.  The main objective of this paper is to examine naturally occurring Banach algebras which are canonically dual spaces and find conditions which imply that the canonical predual is unique. 

A \emph{dual Banach algebra} is a Banach algebra $\mc A$ which is a dual Banach
space such that the product on $\mc A$ is separately weak$^*$-continuous. Important examples include group algebras $\ell^1(G)$ with the canonical predual $c_0(G)$, and their generalisations the group measure algebras, the Fourier algebra and the Fourier-Stieltjes algebra, and the algebra of bounded operators on any reflexive Banach space. In general, the weak$^*$-topology on a dual Banach algebra is not uniquely determined. Indeed, consider any Banach space $E$ which admits two distinct weak$^*$-topologies, such as $\ell^1(\mathbb N)$ for example, and equip it with the zero product to obtain a dual Banach algebra with multiple weak$^*$-topologies.

In \cite{runde}, Runde studied some cohomological properties of dual Banach algebras, although
the idea had been recognised before.  A theory of amenability known as \emph{Connes-amenability}
was developed paralleling that for von Neumann algebras.  It was shown in \cite{Daws} that
any dual Banach algebra is weak$^*$-isomorphic to a weak$^*$-closed subalgebra of $\mc B(E)$,
for some reflexive Banach space $E$.  Furthermore, a characterisation of Connes-amenability
in terms of an injectivity condition was given, again paralleling the von Neumann algebra
theory.  Unfortunately, the non-uniqueness of the weak$^*$-topology means that these concepts
may possibly depend on the choice of the weak$^*$-topology involved.  Our objective in this
paper is to give natural conditions on important examples of dual Banach algebras which ensure
that the weak$^*$-topology is uniquely determined and is induced by the canonical predual.
We show that for a locally compact group $G$, the measure algebra $M(G)$ has a unique
weak$^*$-topology induced by the duality between $C_0(G)$ and $M(G)$ which makes a certain
natural coassociative product weak$^*$-continuous.  That is, $M(G)$ has a unique
weak$^*$-topology as a Hopf algebra. These examples include the group convolution algebras
$\ell^1(G)$ for discrete groups $G$ for which the weak$^*$-topology is not in general unique.
Indeed, with Haydon and Schulmprecht, the first and third author have constructed uncountably
many distinct weak$^*$-topologies on the convolution algebra $\ell^1(\mathbb Z)$ making the
multiplication separately continuous, \cite{DHSW}.  Our second class of examples are the Fourier algebras $A(G)$ for separable compact
groups $G$. In this case we show that there is a unique isometric predual for $A(G)$ making
the Hopf algebra operations continuous.  

Let us outline the ideas involved for the special case of the convolution algebra $\ell^1(G)$ for some countable discrete group $G$.  Any weak$^*$-topology on $\ell^1(G)$ is induced by a closed subspace $E$ of $\ell^\infty(G)=\ell^1(G)'$. Continuity of the coassociative product making $\ell^1(G)$ into a Hopf algebra is equivalent to $E$ being a subalgebra of $\ell^\infty(G)$ with the usual pointwise product.  Some commutative C$^*$-algebra theory allows us to conclude that $E$ must actually be a C$^*$-subalgebra of $\ell^\infty(G)$. Thus the Gelfand transform allows us to identify $E$ with $C_0(K)$ for some locally compact Hausdorff space $K$.  The dual pairing between $E$ and $\ell^1(G)$ allows us to conclude that characters on $E$ must arise from elements of $G$, so there is a natural bijection between $G$ and $K$ inducing a group structure on $K$. We then show that this group structure is compatible with the topology on $K$, that this topology must be discrete, and finally that $E$ must be the subspace $c_0(G)$ of $\ell^\infty(G)$. Thus the weak$^*$-topology induced on $\ell^1(G)$ by $E$ is the canonical one.  This  outline extends to the $M(G)$ case, which we consider in section \ref{measure_case}. 

To consider $A(G)$, the suitable non-commutative analogue of the character space, crucial to the outline above, is the \emph{spectrum} of a type $I$ $C^*$-algebra.  Again, any weak$^*$-topology making the coassociative product continuous is induced by a subalgebra $E$ of $VN(G)$, and if the canonical isomorphism between $A(G)$ and $E$ is isometric, then $E$ is a C$^*$-subalgebra of $VN(G)$, as shown in Theorem \ref{autoc*}. In section \ref{prelim} we examine the general properties of preduals which are also $C^*$-subalgebras in this fashion, giving a bijective correspondence between the representation theory of $E$ and the canonical predual. The main remaining difficulty is to ensure that this bijection is a homemorphism between the spectrums of $E$ and $C^*(G)$. This is pursued in Section~\ref{fourier}.  We need the additional isometric assumption, as our arguments above in the non-isometric case crucially depended upon commutative theory.

In section \ref{Uniqueness}, we consider classes of dual Banach algebras whose weak$^*$-topology is uniquely determined.  In particular, we extend the earlier work of the first author \cite{Daws} to show that a von Neumann algebra has this property, that is a weak$^*$-topology which makes the multiplication separately continuous \emph{automatically} makes the adjoint continuous and so is the usual weak$^*$-topology. 

\smallskip\textbf{Acknowledgements:} Part of this research was undertaken during a visit of the first author to Texas A\&M University in 2007. He would like to thank the faculty of Texas A\&M for their hospitality. We would like to thank Brian Forrest and Nico Spronk for pointing us in the direction of the papers \cite{TB} and \cite{taylor}.

\section{Preliminaries}\label{prelim}

Let $\mc{A}$ be a Banach algebra. For a moment, let us forget the multiplication and just consider $\mc{A}$ as a Banach space. A weak$^*$-topology on $\mc{A}$ is induced by a Banach space $\tilde{E}$ whose dual $\tilde{E}'$ is isomorphic as a Banach space to $\mc A$, via $j:\mc{A}\rightarrow\tilde{E}'$ say.  We do not assume that $j$ is isometric.  Let $E$ be the closed subspace $(j'\circ\kappa_{\tilde{E}})(\tilde{E})$ of $\mc{A}'$ whose dual is canonically identified with $\mc{A}''/E^\perp$, where
$$
E^\perp=\{x\in\mc{A}'':\ip{\mu}{x}=0,\quad\forall\mu\in E\}.
$$
Here, and throughout, we write $\ip{x}{\mu}=\mu(x)$ for $\mu\in \mc{A}'$ and $x\in\mc{A}''$ for the dual pairing between $(\mc{A}',\mc{A}'')$. 

The natural map $\iota_E:\mc{A}\rightarrow E'=\mc{A}''/E^\perp$, arising as the composition of $\kappa_\mc{A}$ and the quotient map onto $\mc{A}''/E^\perp=E'$, is an isomorphism, and the weak$^*$-topologies induced by $\tilde{E}$ and $E$ agree.  Thus it suffices to consider weak$^*$-topologies induced by subspaces $E$ of $\mc{A}'$ for which the natural map $\iota_E:\mc{A}\rightarrow\mc{A}''/E^\perp$ is an isomorphism.  Furthermore, by concretely realising preduals as subspaces of $\mc A'$, we have the following useful fact.  Two preduals $E_0 \subseteq \mc A'$ and $E_1\subseteq \mc A'$ induce distinct weak$^*$-topologies on $\mc A$ if, and only if, $E_0 \not= E_1$ as subspaces.

\begin{lemma}
The map $\iota_E$ is isometric when $j$ is isometric.
\end{lemma}
\begin{proof}
Notice that clearly $\iota_E$ is a contraction.
If $j$ is a surjective isometry, then so is $j':\tilde E''\rightarrow\mc A'$, and so, $j' \kappa_{\tilde E}$ is an isometry onto its range.  Hence, for $a\in\mc A$, we have
\begin{align*} \|\iota_E(a)\| &= \sup\big\{ |\ip{\mu}{a}| : \mu\in E, \|\mu\|\leq1 \big\} \\
&= \sup\big\{ |\ip{j'\kappa_{\tilde E}(x)}{a}| : x\in \tilde E, \|j'\kappa_{\tilde E}(x)\|\leq 1 \big\} \\
&= \sup\big\{ |\ip{j(a)}{x}| : x\in \tilde E, \|x\|\leq 1 \big\} = \|j(a)\| = \|a\|.\qedhere
\end{align*}
\end{proof}

Now let us examine when a weak$^*$-topology makes the multiplication separately continuous. The dual $\mc{A}'$ is naturally an $\mc{A}$-module via the definitions
$$
\ip{a\cdot\mu}{b} = \ip{\mu}{ba}, \quad
\ip{\mu\cdot a}{b} = \ip{\mu}{ab} \qquad (a,b\in\mc A,\mu\in\mc A').
$$
Given a closed subspace $E$ of $\mc{A}'$ such that $\iota_E:\mc{A}\rightarrow\mc{A}''/E^\perp=E'$ is an isomorphism, it is easy to check that the multiplication in $\mc{A}$ is separately weak$^*$-continuous in the weak$^*$-topology induced by $E$ if, and only if, $E$ is a submodule of $\mc{A}'$.  Thus, in considering weak$^*$-topologies on dual Banach algebras $\mc{A}'$ it suffices to consider closed $\mc{A}$-submodules of $\mc{A}'$.  

\begin{definition}
Let $\mc{A}$ be a Banach algebra.  A \emph{predual} for $\mc{A}$ is a closed submodule $E$ of $\mc{A}'$ such that the composition $q\circ\kappa_{\mc{A}}$ is an isomorphism, where $q:\mc{A}''\rightarrow \mc{A}''/E^\perp=E'$ is the quotient map.  Given a predual $E$ for $\mc{A}$, write $\iota_E=q\circ\kappa_{\mc{A}}$ for the natural isomorphism from $\mc{A}$ onto $E'$. We say that $E$ is an \emph{isometric predual} for $\mc{A}$ if $\iota_E$ is an isometric isomorphism.  If $\mc{A}$ has a predual, then $\mc{A}$ is a dual Banach algebra.
\end{definition}

At this point it is worth noting that not all preduals need be isometric. This is most easily seen by resorting to the trivial product on $\ell^1$, however, in the forthcoming paper \cite{DHSW} examples of non-isometric preduals for $\ell^1(\mathbb Z)$ with the usual convolution multiplication will be exhibited.

The two central classes of examples in this paper are the measure algebras $M(G)$ of a locally compact group and the Fourier algebra $A(G)$ of a separable compact group.  Both these algebras are the unique isometric preduals of von Neumann algebras: $M(G)'\cong C_0(G)''$ and $A(G)'\cong VN(G)$, and so any predual for $M(G)$ or $A(G)$ is a submodule of $C_0(G)''$ or $VN(G)$ respectively.  We now consider this situation in more generality.

Let $F$ be a C$^*$-algebra and suppose that the dual $\mc{A}=F'$ is a Banach algebra.  Suppose that $E\subseteq\mc A'$ is a predual for $\mc{A}$ which is in addition a C$^*$-subalgebra of the von Neumann algebra $F''$. Let $\iota_E:\mc A\rightarrow E'$ be the canonical map, which is assumed
to be an isomorphism as $E$ is a predual for $\mc A$.

\begin{proposition}\label{c*iso}
With the notation above, $\iota_E':E''\rightarrow\mc{A}'=F''$ is a $^*$-homomorphism.
In particular $\iota_E$ and $\iota_E'$ are isometries and $E$ is an isometric predual.
\end{proposition}

\begin{proof}
Let $q:\mc{A}''\rightarrow \mc{A}''/E^\perp=E'$ be the quotient map, so that $\iota_E=q\circ \kappa_{\mc{A}}$
and hence $\iota_E' = \kappa_{\mc{A}}' \circ q'$.  Now, if we identify $E''$ with
$E^{\perp\perp} \subseteq F''''=\mc{A}'''$, then $q':E''\rightarrow \mc{A}'''=F''''$
is just the inclusion map, and is the second adjoint of the inclusion map $E\rightarrow \mc{A}'=F''$,
which is a $^*$-homomorphism.  Hence $q'$ is a $*$-homomorphism.  Furthermore,
$\kappa_{\mc{A}}':\mc{A}''' \rightarrow \mc{A}'=F''$ is a $*$-homomorphism, a fact shown
in much more generality by Palmer in \cite{palmer}.
\end{proof}

Thus the isometry $\iota_E:\mc{A}=F'\rightarrow E'$ is a isomorphism between the duals of two C$^*$-algebras. Next we confirm that it respects the order structure.
\begin{proposition}\label{states}
Let $E,F$ and $\mc A$ be as above.  Let $m\in E'$ be a functional,
and let $a=\iota_E^{-1}(m)\in\mc A=F'$.  Then we have the following:
\begin{enumerate}
\item\label{states:one} $a$ is positive on $F$ if ,and only if, $m$ is positive on $E$;
\item\label{states:two} $a$ is a state on $F$ if, and only if, $m$ is a state on $E$;
\item\label{states:three} $a$ is a pure state on $F$ if, and only if, $m$ is a pure state on $E$;
\end{enumerate}
\end{proposition}
\begin{proof}
\begin{enumerate}
\item The functional $a$ is positive on $F$ if, and only if, $\kappa_{F'}(a)$ is positive on $F''$. Since $\iota_E'$ is an injective $^*$-homomorphism from $E''$ onto $F''$, it follows that $\kappa_{F'}(a)$ is positive if, and only if, $(\iota_E''\circ \kappa_{F'})(a)=\kappa_{E'}\circ\iota_E(a)=\kappa_{E'}(m)$ is positive.  

\item This is immediate as $\iota_E'$ is an isometry and a state is a positive linear functional of norm one.

\item A pure state is an extreme point of the state space, and so this follows by linearity and the previous part.\qedhere
\end{enumerate}
\end{proof}

There is also a 1-1-correspondence between $^*$-representations of $E$ and $F$.    Following the notation of \cite[Section~III.2]{tak}, given a representation $\pi:
E\rightarrow\mc B(H)$ of the C$^*$-algebra $E$ on some Hilbert space, there is a canonical representation $\tilde\pi:E''\rightarrow\mc B(H)$.  This is the unique normal representation of $E''$ satisfying $\tilde\pi\circ\kappa_E=\pi$.  Take $u,v\in H$, and define a functional $\omega(\pi;u,v)$ on $E$
by
\[ \ip{\omega(\pi;u,v)}{x} = (\pi(x)(u)|v) \qquad (x\in E). \]
Here, and throughout, we use the notation $(\cdot|\cdot)$ for the inner product on a Hilbert space, to avoid confusion with the dual pairing used elsewhere.  Then $\tilde\pi$ is defined by
\[ (\tilde\pi(x)(u)|v) = \ip{x}{\omega(\pi;u,v)} \qquad (x\in E'',
u,v\in H). \]
We can then define a map $\phi:F\rightarrow\mc B(H)$ by
$\phi = \tilde\pi\circ(\iota_E^{-1})'\circ\kappa_F$.  This is certainly a representation, since $\iota_E'$ is a $^*$-isomorphism.  We then apply the process above again to form the canonical extension $\tilde\phi$ to $F''$, which is the unique normal representation of $F''$ satisfying $\tilde\phi\circ\kappa_F=\phi$.  This uniqueness ensures that the restriction of $\tilde\phi$ to $E$ is equal to $\phi$.  Let us state this formally for later reference.
\begin{proposition}\label{reps}
With notation as above, $\phi$ is a $*$-representation,
and $\tilde\phi$ restricted to $E$ is equal to $\pi$.
\end{proposition}

Next we consider preduals $E\subset\mc{A}'$ which are a subalgebras of the von Neumann algebra $F''$.  In the situations of interest later, this will occur precisely when the natural co-associative products we consider on $\mc{A}$ are weak$^*$-continuous on the topology induced by $E$.  Under the additional assumption that $E$ is an isometric predual, the next result enables us to reduce to the case when $E$ is a C$^*$-subalgebra and hence use the machinery developed above.  

\begin{theorem}\label{autoc*}
Suppose that $F$ is a C$^*$-algebra, and that $\mc{A}=F'$ is a Banach algebra whose multiplication is separately $\sigma(F',F)$-continuous.  If $E\subset\mc{A}'$ is an isometric predual for $\mc{A}$ making the multiplication separately continuous, and $E$ is also a subalgebra of $F''=\mc{A}'$, then $E$ is a C$^*$-subalgebra of $F''$.
\end{theorem}
\begin{proof}
Let $j$ denote the inclusion map from $E$ into $F''$ so that $j=\iota_E'\circ\kappa_E$, where $\iota_E$ is the canonical isomorphism from $\mc{A}$ to $E'$ which is an isometry by hypothesis.  Since $\iota_E$ factorises as $q\circ\kappa_{\mc{A}}$, where $q:F'''\rightarrow F'''/E^\perp$ it follows that $\iota_E'$ is an algebra homomorphism, just as in Proposition \ref{c*iso}. Thus $\alpha=j''\circ(\iota_E')^{-1}$ is an isometric algebra homomorphism from $F''$ into $F''''$.  Such maps are necessarily $^*$-preserving.  Indeed, let $1_{F''}$ denote the unit of $F''$, so that $p=\alpha(1_{F''})$ is an idempotent of norm one in $F''''$ and therefore is a projection. Let $B$ be the C$^*$-subalgebra of $F''''$ generated by $\alpha(F'')$, which by the preceding calculation, has identity $p=\alpha(1_{F''})$.  Given a unitary $u\in F''$, it follows that
$$
\alpha(u)\alpha(u^*)=\alpha(u^*)\alpha(u)=p.
$$
Since $\|\alpha(u)\|=\|\alpha(u^*)\|=1$, it follows that $\alpha(u)$ is a unitary element in $B$ with $\alpha(u)^*=\alpha(u^*)$.  As $F''$ is spanned by its unitaries, $\alpha$ is $^*$-preserving as claimed.

The pre-adjoint $\beta:F'''\rightarrow F'$ of the isometric $^*$-homomorphism $\alpha$ is also a $^*$-preserving map.  Given any $\phi\in F'''$, it follows that $\phi|_{j(E)}=0$ if, and only if, $\beta(\phi)=0$ since
$$
\ip{\beta(\phi)}{x}=\ip{\phi}{\alpha(x)}=\ip{\phi}{(j''\circ(\iota_E')^{-1})(x)}=0,\quad x\in F''.
$$
Thus, if $\phi|_{j(E)}=0$, then $\beta(\phi)=0$ and so $\beta(\phi^*)=0$ when $\phi^*|_{j(E)}=0$.  Therefore 
$$
\phi|_{j(E)}=0\Longrightarrow\phi|_{j(E)^*}=0,\quad \phi\in F''',
$$
so the Hahn-Banach theorem implies that $j(E)=j(E)^*$, as required.
\end{proof}

A careful examination of the proof above reveals that we do not need $\mathcal A$ to be a Banach algebra in Theorem \ref{autoc*}.  Thus any isometric predual for the dual of a C$^*$-algebra which is also a subalgebra of the bidual is automatically a C$^*$-subalgebra of this bidual.

\section{Group measure algebras}\label{measure_case}

Let $G$ be a locally compact group, and form the Banach space $C_0(G)$ of
continuous functions on $G$ which vanish at infinity.  The dual is
$M(G)$, the space of regular finite measures on $G$, which becomes a Banach
algebra under the convolution product, defined by
\[ \ip{\mu\lambda}{f} = \int_{G\times G} f(st) \, d\mu(s) \, d\lambda(t)
\qquad (\mu,\lambda\in M(G), f\in C_0(G)). \]
For each $s\in G$, let $\delta_s$ be the point mass measure at $s$, so $\delta_s$
is the character on $C_0(G)$ formed by evaluation at $s$.  Then
$\delta_s \delta_t = \delta_{st}$ for each $s,t\in G$, and the family
$\{ \delta_s : s\in G \}$ has weak$^*$-dense linear span in $M(G)$.

An extension of Wendel's Theorem (see \cite[Theorem~3.3.40]{Dales}) due to Johnson \cite{johnson} shows that $M(G)$ is an isometric invariant for the locally compact group $G$.  However, for example, $\ell^1(C_4)$ and $\ell^1(C_2\times C_2)$ are isomorphic Banach algebras, but are not isometric.  To more fully capture the group structure,
we introduce a coassociative product, now a common concept from the
theory of quantum groups.  This idea was developed in \cite{ER} which we
shall turn to later when considering the Fourier algebra, but for $M(G)$, we can sketch the theory in an elementary way.

For a locally compact group $G$, we define
\[ \Gamma = \Gamma_G : M(G)\rightarrow M(G\times G);\quad
\delta_s \mapsto \delta_{(s,s)} \qquad (s\in G). \]
To be more precise,
\[ \ip{\Gamma(\mu)}{f} = \int_G f(s,s) \, d\mu(s)
\qquad (\mu\in M(G), f\in C_0(G\times G)). \]
Notice that for $f\in C_0(G\times G)$ and $\mu,\nu\in M(G)$,
\begin{align*} \ip{\Gamma(\mu\nu)}{f} &= \int_G f(s,s) \, d(\mu\nu)(s)
= \int_G \int_G f(st,st) \, d\mu(s) \, d\nu(t), \\
&= \int_{G\times G} \int_{G\times G} f(su,tv) \, d\Gamma(\mu)(s,t) \, d\Gamma(\nu)(u,v)
= \ip{\Gamma(\mu)\Gamma(\nu)}{f}. \end{align*}
Hence $\Gamma$ is a Banach algebra homomorphism for the convolution
products on $M(G)$ and $M(G\times G)$.

It seems natural to insist that this map be weak$^*$-continuous,
that is, there should exist a map $\Gamma_*:C_0(G\times G)
\rightarrow C_0(G)$ with $\Gamma_*'=\Gamma$.  Clearly, for
$f\in C_0(G\times G)$, we have that $\Gamma_*(f)=g\in C_0(G)$ where
$g(s)=f(s,s)$ for $s\in G$.  We identify $C_0(G)\otimes C_0(G)$
with a dense subspace of $C_0(G\times G)$ in the natural way,
and then for $f,g\in C_0(G)$, we see that
\[ \Gamma_*(f\otimes g)(s) = (f\otimes g)(s,s) = f(s)g(s)
\qquad (s\in G). \]
Hence $\Gamma_*$ induces the usual multiplication on $C_0(G)$. It is no surprise that with this additional structure, $M(G)$ fully captures the group $G$.

\begin{proposition}\label{meas_coass}
Let $G$ and $H$ be locally compact groups, and let $\theta:M(G)
\rightarrow M(H)$ be a Banach algebra isomorphism.  Suppose furthermore
that $\theta$ intertwines the Hopf-algebra products, that is,
$(\theta\otimes\theta)\Gamma_G = \Gamma_H\theta$.  Then there is
a group isomorphism $\alpha:G\rightarrow H$ such that $\theta=\alpha^*$,
that is,
\[ \ip{\theta(\mu)}{f} = \int_G f(\alpha(s)) \ d\mu(s)
\qquad (\mu\in M(G), f\in C_0(H)). \]
In particular, $\theta$ is automatically weak$^*$-continuous.
\end{proposition}
\begin{proof}
Let $f,g\in C_0(H)''$, and let $\mu\in M(G)$, so we see that
\begin{align*} \ip{\theta'(fg)}{\mu} &= \ip{f\otimes g}{\Gamma_H\theta(\mu)}
= \ip{f\otimes g}{(\theta\otimes\theta)\Gamma_G(\mu)} \\
&= \ip{\theta'(f) \otimes \theta'(g)}{\Gamma_G(\mu)}
= \ip{\theta'(f) \theta'(g)}{\mu}.
\end{align*}
Hence $\theta':C_0(H)''\rightarrow C_0(G)''$ is an algebra homomorphism.
As $\theta'$ is also an isomorphism, we conclude that $\theta'$ is
an isometry.  So $\theta:M(G)\rightarrow M(G)$ is an isometric isomorphism,
and so by Johnson's Theorem, \cite{johnson}, there exists a continuous character
$\chi:G\rightarrow(0,\infty)$ and a topological group isomorphism
$\alpha:G\rightarrow H$, such that
\[ \ip{\theta(\mu)}{f} = \int_G f(\alpha(s)) \chi(s) \, d\mu(s)
\qquad (\mu\in M(G), f\in C_0(H)). \]
In particular, $\theta(\delta_s) = \chi(s) \delta_{\alpha(s)}$ for each $s\in G$.
For $f,g\in C_0(H)$, we hence have that
\[ \ip{\theta'(fg)}{\delta_s} = \chi(s) f(\alpha(s)) g(\alpha(s))
= \ip{\theta'(f)}{\delta_s}\ip{\theta'(g)}{\delta_s} = \chi(s)^2
f(\alpha(s)) g(\alpha(s)). \]
Consequently $\chi(s)=\chi(s)^2$ for each $s\in G$.  Thus $\chi\equiv 1$, which
completes the proof.
\end{proof}

Now suppose that $E\subseteq M(G)'$ is a predual for $M(G)$.
We need a slight digression into the theory tensor products of Banach spaces,
for which we refer the reader to \cite{Ryan}.  We form the
\emph{injective tensor product} of $E$ with itself, $E \inten E$.
\begin{proposition}
With notation as above, $(E\inten E)'$ is naturally isomorphic $M(G\times G)$.
\end{proposition}
\begin{proof}
The dual of $E\inten E$ can be identified with the space of \emph{integral}
operators $E\rightarrow E'$, written $\mc I(E,E')$, by the dual pairing
\[ \ip{T}{x\otimes y} = \ip{T(x)}{y} \qquad (x\otimes y\in E\inten E). \]
We identify $C_0(G) \inten C_0(G)$ with $C_0(G\times G)$, and so identify
$\mc I(C_0(G),M(G))$ with $M(G\times G)$ by, for $\mu\in M(G\times G)$,
we define $T_\mu\in\mc I(C_0(G),M(G))$ by
\[ \ip{T_\mu(f)}{g} = \ip{\mu}{f\otimes g} \qquad (f,g\in C_0(G)). \]

Given $\mu\in M(G\times G)$, as $T_\mu$ is weakly-compact, $T_\mu$ takes
$C_0(G)''$ into $\kappa_{M(G)}(M(G))$, and so there exists $T_0\in
\mc B(C_0(G)'',M(G))$ such that $\kappa_{M(G)} T_0 = T_\mu''$.  Then
$T_0 = \kappa_{C_0(G)}' \kappa_{M(G)} T_0 = \kappa_{C_0(G)}'T_\mu''$ is
integral, as the integral operators form an operator ideal.  Define
$S_\mu\in\mc B(E,E')$ by $S_\mu = \iota_E T_0 \iota_E' \kappa_E$, so that
$S_\mu$ is integral, hence $S_\mu \in (E\inten E)'$.  Let $x\in E\subseteq M(G)'$,
so that $T_\mu''(x) = \kappa_{M(G)}(\mu)$ say, for some $\mu\in M(G)$.  We
can verify that $\iota_E'\kappa_E:E\rightarrow M(G)'$ is just the inclusion
map, and so, for $y\in E$,
\begin{align*} \ip{S_\mu}{x\otimes y} &= \ip{S_\mu(x)}{y} = \ip{S_\mu'\kappa_E(y)}{x}
= \ip{\kappa_E' \iota_E'' T_0' \iota_E' \kappa_E(y)}{x} \\
&= \ip{\kappa_E' \iota_E'' T_\mu''' \kappa_{C_0(G)}'' \iota_E' \kappa_E(y)}{x}
= \ip{\iota_E \kappa_{C_0(G)}' T_\mu'' \iota_E' \kappa_E(x)}{y} \\
&= \ip{\iota_E \kappa_{C_0(G)}'\kappa_{M(G)}(\mu)}{y}
= \ip{\iota_E(\mu)}{y} = \ip{T_\mu''(x)}{y}.
\end{align*}
We have hence defined a map $\phi:M(G\times G) \rightarrow (E\inten E)'$ by
$\phi(\mu)=S_\mu$.

Similarly, define $\psi:(E\inten E)'=\mc I(E,E') \rightarrow
\mc I(C_0(G),M(G))=M(G\times G)$ by
\[ \psi(S) = \kappa_{C_0(G)}' (\iota_E^{-1})'' S'' (\iota_E^{-1})' \kappa_{C_0(G)}
\qquad ( S\in\mc I(E,E') ). \]
Again, this is well-defined, as if $S\in\mc I(E,E')$ then
$S'\in\mc I(E'',E')$ and so $\psi(S)\in\mc I(C_0(G),M(G))$.  Let $\mu\in M(G\times G)$,
so that
\begin{align*}
\iota_E^{-1} S_\mu' (\iota_E^{-1})' \kappa_{C_0(G)} &=
\iota_E^{-1} \kappa_E' \iota_E'' T_0' \iota_E' (\iota_E^{-1})' \kappa_{C_0(G)}
= \iota_E^{-1} \kappa_E' \big( \iota_E T_\mu' \kappa_{C_0(G)}\big)'' \kappa_{C_0(G)} \\
&= \iota_E^{-1} \kappa_E' \kappa_{E'} \big( \iota_E T_\mu' \kappa_{C_0(G)}\big)
= \iota_E^{-1} \iota_E T_\mu' \kappa_{C_0(G)} = T_\mu'\kappa_{C_0(G)},
\end{align*}
and so
\[ \psi(S_\mu) = \kappa_{C_0(G)}' T_\mu'' \kappa_{C_0(G)}
= \kappa_{C_0(G)}' \kappa_{M(G)} T_\mu = T_\mu. \]
Hence $\psi = \phi^{-1}$, and we see that $E\inten E$ is a predual for $M(G\times G)$.
\end{proof}

The key idea is that given a predual $E$ for $M(G)$, we have a
natural way of forming a predual for $M(G\times G)$, namely
taking some completion of the tensor product $E\otimes E$.  It hence
makes sense to ask when $\Gamma$ is weak$^*$-continuous for the predual $E$.  Let us call a predual making $\Gamma$ weak$^*$-continuous a \emph{Hopf algebra predual of $M(G)$}, and the induced weak$^*$-topology for $M(G)$ a \emph{Hopf-algebra weak$^*$-topology}. Firstly let us characterise these preduals algebraically.
\begin{lemma}
Let $E\subset M(G)'=C_0(G)''$ be a predual for $M(G)$.  Then $\Gamma$ is weak$^*$-continuous if, and only if, $E$ is a subalgebra of the von Neumann algebra $C_0(G)''$.
\end{lemma}
\begin{proof}
If $\Gamma$ is weak$^*$-continuous,
then there exists $\Gamma_*:E\inten E\rightarrow E$ such that $\Gamma_*'=\Gamma$.
Let $\mu\in M(G)$, so
\[ \ip{\Gamma(\mu)}{f\otimes g} = \int_G f(s) g(s) \, d\mu(s)
\qquad (f,g\in C_0(G)), \]
so we see that, in the notation of the previous proof, $T_{\Gamma(\mu)}(f) = f\mu$
where $f\mu\in M(G)$ is the measure given by $\ip{f\mu}{g} = \ip{\mu}{fg}$
for $g\in C_0(G)$.  Let $x,y\in E$, so we can find bounded nets $(f_\alpha)$
and $(g_\alpha)$ in $C_0(G)$ such that
\[ \ip{x}{\lambda} = \lim_\alpha \ip{\lambda}{f_\alpha}, \quad
\ip{y}{\lambda} = \lim_\alpha \ip{\lambda}{g_\alpha} \qquad (\lambda\in M(G)). \]
We then see that
\begin{align*}
\ip{\Gamma_*(x\otimes y)}{\mu} &= \ip{x\otimes y}{\Gamma(\mu)}
= \ip{x}{T_{\Gamma(\mu)}'(y)} = \lim_\beta \ip{T_{\Gamma(\mu)}'(y)}{f_\beta} \\
&= \lim_\beta \lim_\alpha \ip{T_{\Gamma(\mu)}(f_\beta)}{g_\alpha}
= \lim_\beta \lim_\alpha \ip{\mu}{f_\beta g_\alpha}
= \ip{xy}{\mu},
\end{align*}
by the definition of the Arens products (see after \cite[Theorem~2.6.15]{Dales}
for example).  Consequently, we see that $\Gamma_*$ maps into $E$ only
when the multiplication on $C_0(G)''$ restricts to $E$, that is, $E$ is
a subalgebra of $C_0(G)''$.
\end{proof}

In fact Hopf-algebra preduals for $M(G)$ are automatically isometric.

\begin{lemma}\label{cstar_bidual_is_cstar}
Let $K$ be a locally compact Hausdorff space, let $L$ be a compact Hausdorff
space, and let $\mc A$ be a closed subalgebra of $C_0(K)$ such that $\mc A''$
is Banach algebra isomorphic to $C(L)$.  Then $\mc A$ is a C$^*$-subalgebra
of $C_0(K)$.
\end{lemma}
\begin{proof}
Let $\mc B$ be the C$^*$-subalgebra of $C_0(K)''$ generated by $\mc A''$.
Thus $\mc B$ is unital, and so we identify $\mc B$ with $C(S)$ for some compact
Hausdorff space $S$.  Let $\theta:C(L) \rightarrow \mc A''$ be an algebra isomorphism,
so we can regard $\theta$ as an algebra homomorphism
$C(L)\rightarrow C(S)$.  Hence there exists a continuous map $\alpha:S\rightarrow L$ such that
\[ \theta(f)(s) = f(\alpha(s)) \qquad (f\in C(L), s\in S). \]
As $\mc A''$ separates the points of $S$, we see that $\alpha$ is injective.
As $S$ is compact, $\alpha$ is a homeomorphism onto its range.  As $\theta$ is
an isomorphism, $\alpha$ must be a homeomorphism between $S$ and $L$.
We hence conclude that $\mc A'' = \mc B = C(S)$.

Suppose that $\mc A$ is not a C$^*$-subalgebra of $C_0(K)$, so there exists
$a_0\in\mc A$ such that $a_0^*\not\in\mc A$.  Hence there exists $\lambda\in C_0(K)'$
such that $\ip{\lambda}{a_0^*}=1$ and $\ip{\lambda}{a}=0$
for all $a\in\mc A$.  As $\lambda$ is a normal functional on $C_0(K)''$,
and so we can regard $\lambda$ as a functional on $C(S)$.  As $\lambda$ kills
$\mc A$, $\lambda$ also kills $\mc A''=C(S)$.  Hence also $\lambda^*$ kills
$C(S)$.  However, we then see that
\[ 1 = \ip{\lambda}{a_0^*} = \ip{\lambda^*}{a_0} = 0, \]
a contradiction, as $a_0 \in \mc A\subseteq C(S)$.  So $\mc A$ is a C$^*$-subalgebra of $C_0(K)$.
\end{proof}

\begin{corollary}\label{HopfM2}
Any Hopf-algebra predual for $M(G)$ is automatically an isometric predual
and is a C$^*$-subalgebra of $C_0(G)''$.  
\end{corollary}
\begin{proof}
Let $E$ be a Hopf algebra predual of $M(G)$, so $E$ is a subalgebra of $C_0(G)''$.
Let $\iota_E:M(G)\rightarrow E'$ be the canonical isomorphism.  Using a similar argument
to that used in the proof of Proposition~\ref{c*iso}, we can show that $\iota_E'$ is
an algebra homomorphism (see also \cite{palmer}).
It follows that $E''$ and $C_0(G)''$ are isomorphic as Banach algebras, and so
by Lemma \ref{cstar_bidual_is_cstar}, $E$ is a C$^*$-subalgebra of $C_0(G)''$,
and hence is an isometric predual by Proposition \ref{c*iso}.
\end{proof}

Our main result in this section is that there is a unique Hopf algebra weak$^*$-topology on $M(G)$.

\begin{theorem}\label{measure.hopf}
The canonical predual $C_0(G)$ is the unique Hopf algebra predual of $M(G)$.
\end{theorem}
\begin{proof}
Let $E\subset C_0(G)''$ be a Hopf-algebra predual for $M(G)$.  By Corollary~\ref{HopfM2}, $E$ is a C$^*$-subalgebra of $C_0(G)''$. Let $\iota_E:M(G)\rightarrow E'$ be the canonical isomorphism. Let $K$ denote the character space of $E$ and let $\mathcal G:E\rightarrow C_0(K)$ be the Gelfand transform.  Identify the characters on $C_0(G)$ with $G$.  By Lemma \ref{states}, $\iota_E^{-1}$ restricts to a bijection $\theta:G\rightarrow K$ and we use this bijection to induce a group structure on $K$. Notice that for $f\in C_0(K)$ and $k\in K$,
\[ f(k) = \ip{\delta_k}{f} = \ip{\mc G^{-1}(f)}{\delta_{\theta^{-1}(k)}}. \]

We now claim that $K$ is a \emph{semitopological semigroup}, i.e. the product induced by $\theta$ is separately continuous.  Let $(k_\alpha)$ be a net in $K$ converging to $k$, and let $l\in K$.
Then, for $f\in C_0(K)$, using the fact that $E\subseteq M(G)'$ is an
$M(G)$-bimodule, we have
\begin{align*}
\lim_\alpha \ip{\delta_{k_\alpha l}}{f}
&= \lim_\alpha \ip{\mc G^{-1}(f)}{\delta_{\theta^{-1}(k_\alpha)} \delta_{\theta^{-1}(l)}}
= \lim_\alpha \ip{\delta_{\theta^{-1}(l)} \cdot \mc G^{-1}(f)}{\delta_{\theta^{-1}(k_\alpha)}} \\
&= \lim_\alpha \ip{\mc G\big( \delta_{\theta^{-1}(l)} \cdot \mc G^{-1}(f) \big)}{\delta_{k_\alpha}}
= \ip{\mc G\big( \delta_{\theta^{-1}(l)} \cdot \mc G^{-1}(f) \big)}{\delta_k} \\
&= \ip{\mc G^{-1}(f)}{\delta_{\theta^{-1}(k) \theta^{-1}(l)}} = \ip{\delta_{kl}}{f}.
\end{align*}
Hence $k_\alpha l \rightarrow kl$. Analogously, $l k_\alpha\rightarrow lk$, which establishes the claim.

Ellis's Theorem, see \cite{ellis}, says that any locally compact semitopological
semigroup which is algebraically a group is in fact a topological group, that is,
the product is jointly continuous and the inverse is continuous.  In particular
$K$, equipped with the product induced by $\theta$ is a locally compact
topological group.  Now we show that $\theta$ is also a homeomorphism from
$G$ to $K$ and so is a topological group isomorphism.

Define $T=(\mc G')^{-1} \iota_E:M(G)\rightarrow M(K)$. Since $\iota_E$ is an isometric isomorphism, $T$ is an isometric isomorphism of Banach spaces.  Let $(\mu_\alpha)$
be a net in $M(G)$ which converges to $\mu$ in the weak$^*$-topology
induced by $E$.  Then, for $f\in C_0(K)$,
\[ \lim_\alpha \ip{T(\mu_\alpha)}{f} = \lim_\alpha \ip{\mc G^{-1}(f)}{\mu_\alpha}
= \ip{\mc G^{-1}(f)}{\mu} = \ip{T(\mu)}{f}. \]
Thus $T$ is weak$^*$-continuous.
For $f\in C_0(K)$ and $k\in K$, we have that
\[ \ip{\delta_k}{f} = \ip{\mc G^{-1}(f)}{\delta_{\theta^{-1}(k)}}
= \ip{T( \delta_{\theta^{-1}(k)} )}{f}, \]
so that $\delta_{\theta(s)} = T( \delta_s )$ for $s\in G$.
By weak$^*$-continuity and density, we conclude that $T$ is also an algebra
homomorphism.  By Johnson's result in \cite{johnson}
(which itself follows from Wendel's Theorem, see \cite[Theorem~3.3.40]{Dales})
there exists a continuous character $\chi$ on $G$ and a topological
group isomorphism $\psi:G\rightarrow K$, such that
\[ \ip{T(\mu)}{f} = \int_G f(\psi(t)) \chi(t) \, d\mu(t)
\qquad (\mu\in M(G), f\in C_0(K)). \]
Hence, for $f\in C_0(K)$ and $k\in K$, we have that
\begin{align*}
f(k) = \int_G f(\psi(t)) \chi(t) \, d\delta_{\theta^{-1}(k)}(t)
= f(\psi\theta^{-1}(k)) \chi(\theta^{-1}(k)), \end{align*}
so we conclude that $\chi=1$ identically, and that $\psi = \theta$.

Let $x\in E$, and define $f:G\rightarrow\mathbb C$ by
\[ f(t) = \ip{x}{\delta_t} \qquad (t\in G). \]
Then $f(t) = \mc G(x)(\theta(t))$, and so $f\in C_0(G)$.
For $\mu\in M(G)$, we see that
\[ \ip{x}{\mu} = \ip{T(\mu)}{\mc G(x)}
= \int_G \mc G(x)(\theta(t)) \, d\mu(t) = \ip{\mu}{f}. \]
Thus $x = f \in C_0(G)$, and we conclude that $E=C_0(G)$.
\end{proof}

\section{Fourier algebras}\label{fourier}

In this section we shall study the `Fourier transform' of the $\ell^1$
group algebra for a discrete group, namely the Fourier algebra of a
compact group.  Let $G$ be a locally compact group, and let $\hat G$
be the collection of equivalence classes of irreducible unitary
representations of $G$.  Recall that every unitary representaion
$\pi:G\rightarrow\mc B(H)$ extends to a $*$-representation $\pi:L^1(G)
\rightarrow\mc B(H)$, and that $C^*(G)$, the group C$^*$-algebra of $G$,
is the completion of $L^1(G)$ with respect to the norm
\[ \|f\|_{C^*(G)} = \sup\{ \|\pi(f)\| : \pi\in\hat G \}
\qquad (f\in L^1(G)). \]
The dual space of $C^*(G)$ is $B(G)$, the Fourier-Stieltjes algebra
of $G$, which can be identified as the space of coefficient functions
\[ G\rightarrow\mathbb C; \ g\mapsto (\pi(g)u|v)
\qquad (g\in G), \]
where $\pi:G\rightarrow\mc B(H)$ is a unitary representation, and
$u,v\in H$.  Then $B(G)$ is a subalgebra of $C(G)$ and the product is
given by tensoring unitary representations.

In \cite{Eymard}, Eymard defined $A(G)$, the Fourier algebra of $G$, to be
the closure in $B(G)$ of those functions with compact support.
Alternatively, consider the left-regular
representation $\lambda:G\rightarrow\mc B(L^2(G))$,
\[ \lambda(s)(f)=g, \quad g(t)=f(s^{-1}t)
\qquad (f\in L^2(G), s,t\in G). \]
Then $A(G)$ is the space of coefficient functions associated to
$\lambda$.  The dual of $A(G)$ may be identified with $VN(G)$,
the group von Neumann algebra, which is the von Neumann algebra in
$\mc B(L^2(G))$ generated by $\{ \lambda(s) : s\in G\}$.  The norm closure
of $\lambda(L^1(G))$ in $\mc B(L^2(G))$ is the reduced group $C^*$-algebra
$C^*_r(G)$.

When $G$ is a compact group,
\[ C^*_r(G) = C^*(G), \quad A(G) = B(G) = C^*(G)', \quad
VN(G) = A(G)' = C^*(G)''. \]
In this section, we shall investigate weak$^*$-topologies on $A(G)$.

In contrast to the $\ell^1$ case, and perhaps surprisingly,
it is possible for $A(G)$ to be a dual Banach space when $G$ is
not compact.  For example, in \cite{taylor}, K.~Taylor shows that
when $G$ is separable and has the [AR] property, then $A(G)$ has
a predual (which can be taken to be a C$^*$-algebra).  A group $G$ has
the [AR] property if and only if $VN(G)$ is atomic, and for example,
it is shown in \cite{TB} that the ``$ax+b$'' group is a non-compact
group with the [AR] property.

The following argument was suggested to us by Brian Forrest.  When $G$ is an amenable
group and $A(G)$ is a dual Banach algebra, then by Leptin's theorem $A(G)$
has a bounded approximate identity, and so by taking weak$^*$-limits,
$A(G)$ is unital, and hence $G$ is compact. We do not know of an example
of a non-amenable $G$ for which $A(G)$ is a dual Banach algebra.
Henceforth, we consider the case of compact $G$.

The natural co-multiplication on $A(G)$ is simply the pre-adjoint
of the multiplication on $VN(G)$, say $m:VN(G)\proten VN(G)
\rightarrow VN(G)$.  The naive (although natural) norm to
consider on $A(G) \otimes A(G)$ is the predual of the von
Neumann tensor norm on $VN(G) \overline\otimes VN(G) = 
VN(G\times G)$.  In \cite[Section~5]{Quigg}, Quigg considers
exactly this problem.  He shows that $m$ has a bounded pre-adjoint,
with respect to this norm, only when $G$ is compact, and the
irreducible representations of $G$ have uniformly bounded
dimension (which of course excludes even $G = SU(2)$).
In \cite{ER}, Effros and Ruan work out the details in a manner which handles all compact groups by introducing the \emph{extended Haagerup tensor product} on $A(G)\otimes A(G)$.

Let us briefly develop this theory.  For a von Neumann algebra $\mc M$, the
\emph{normal Haagerup tensor product} $\mc M \norhaa \mc M$
is such that the multiplication $\mc M\otimes\mc M\rightarrow\mc M$ extends to
a weak$^*$-continuous (complete) contraction $m:\mc M \norhaa \mc M\rightarrow\mc M$.
The predual of $\mc M \norhaa \mc M$ is $\mc M_* \exthaa \mc M_*$, the
\emph{extended Haagerup tensor product} of $\mc M_*$ with itself.  Hence, in particular,
we get the preadjoint of the multiplication, $m_*:A(G) \rightarrow A(G) \exthaa A(G)$.
In \cite{ER}, it is shown that $m_*$ is an algebra homomorphism.

Let $G$ be a compact group, and let $E\subseteq VN(G)$ be a predual for $A(G)$.
Then $E$ becomes an operator space by restricting the natural operator
space structure on $VN(G)$, and so we can form $E \otimes^h E$,
the \emph{Haagerup tensor product} of $E$ with itself.
Then \cite[Theorem~5.3]{ER} shows that $(E \otimes^h E)' = E' \exthaa E'$.
It is easily checked that $\iota_E:A(G)\rightarrow E'$ is completely contractive, and so
$\iota_E^{-1}$ is completely bounded.  The discussion before \cite[Lemma~5.2]{ER}
constructs a linear map $\iota_E^{-1} \exthaa \iota_E^{-1}:E' \exthaa E'
\rightarrow A(G) \exthaa A(G)$ which is easily seen to be a complete isomorphism.
Thus we can naturally identify $E \otimes^h E$ as a predual for $A(G) \exthaa A(G)$.

We have that $A(G) \exthaa A(G)$ can be identified with the space of normal functionals
on $VN(G)\otimes^h VN(G)$ (see Page~143 in \cite{ER}, and noting
the misprint there).  Hence the inclusion map $VN(G) \otimes^h VN(G) \rightarrow
VN(G)\norhaa VN(G)$ is a complete isometry.  As the Haagerup tensor product
is injective, we see that the inclusion map $E \otimes^h E \rightarrow VN(G) \otimes^h VN(G)$
is a complete isometry, and so the natural embedding
$E \otimes^h E \rightarrow VN(G) \norhaa VN(G)$ is a complete isometry.  

\begin{proposition}\label{fourier_hopf=>alg}
With the notation above, we have that $E$ makes the product map
$m_*:A(G) \rightarrow A(G) \exthaa A(G)$ weak$^*$-continuous if and
only if $E$ is a (possibly not self-adjoint) subalgebra of $VN(G)$.
\end{proposition}
\begin{proof}
Suppose that $m_*:A(G) \rightarrow A(G) \exthaa A(G)$ is weak$^*$-continuous,
so there exists a map $m_{**}:E \otimes^h E \rightarrow E$ making the 
diagram below commute.
\[ \xymatrix{ E' \ar[rrr]^{m_{**}'} \ar[d]^{\iota_E^{-1}} &&& (E \otimes^h E)' = E' \exthaa E'
\ar[d]^{\iota_E^{-1}\exthaa\iota_E^{-1}} \\
A(G) \ar[rrr]^{m_*} &&& A(G) \exthaa A(G) } \]
Dualising, the diagram below also commutes.
\[ \xymatrix{ E'' &&& (E\otimes^h E)'' = E''\norhaa E'' \ar[lll]_{m_{**}''} \\
VN(G) \ar[u]_{(\iota_E')^{-1}} &&& VN(G) \norhaa VN(G)
\ar[lll]_m \ar[u]_{(\iota_E')^{-1} \otimes (\iota_E')^{-1}} } \]
So let $x,y\in E$, and let $j = \iota_E' \circ \kappa_E:E\rightarrow VN(G)$
be the inclusion map.  Then $(\iota_E')^{-1}j = \kappa_E$ and so
\begin{align*} \kappa_E m_{**}(x\otimes y) &= m_{**}''( \kappa_E(x) \otimes \kappa_E(y) )
= m_{**}''( (\iota_E')^{-1} \otimes (\iota_E')^{-1} )(j(x)\otimes j(y)) \\
&= (\iota_E')^{-1} m (j(x)\otimes j(y)). \end{align*}
Hence $jm_{**}(x\otimes y) = m(j(x)\otimes j(y))$, and so $m_{**}$ is just
the multiplication map induced by $VN(G)$.  We conclude that $E$ is thus
a subalgebra (but maybe not self-adjoint).

The converse is now simply a matter of reversing the argument.
\end{proof}

Just as with the measure algebra, the Hopf algebra structure on $A(G)$ fully captures the group.
\begin{proposition}
Let $G$ and $H$ be compact groups, and let $\theta:A(G)\rightarrow A(H)$
be a Banach algebra isomorphism.  Suppose furthermore that $\theta$
intertwines the coassociative products on $A(G)$ and $A(H)$.  Then
$\theta$ is an isometry, and there exists a bicontinuous group isomorphism
$\phi:H\rightarrow G$ such that
\[ \theta(a)(s) = a(\phi(s)) \qquad (a\in A(G), s\in H). \]
\end{proposition}
\begin{proof}
As in the proof of Proposition~\ref{meas_coass}, as $\theta$ intertwines
the coassociative products on $A(G)$ and $A(H)$, we see that
$\theta':VN(H)\rightarrow VN(G)$ is an algebra homomorphism.  Let
$s\in H$, and consider the map
\[ A(G)\rightarrow\mathbb C; \quad a\mapsto \theta(a)(s) \qquad (a\in A(G)). \]
This is a character on $A(G)$, and so is evaluation at some point $\phi(s)\in G$,
say.  Hence $\theta(a)(s) = a(\phi(s))$ for $a\in A(G)$ and $s\in H$.
Let $\lambda_G:G\rightarrow VN(G)$ and $\lambda_H:H\rightarrow VN(H)$ be the left-regular
representations, so that
\[ \theta(a)(s) = \ip{\lambda_H(s)}{\theta(a)} = \ip{\theta'\lambda_H(s)}{a}
= a(\phi(s)) = \ip{\lambda_G(\phi(s))}{a} \qquad (a\in A(G), s\in H). \]
Hence $\theta'\lambda_H(s) = \lambda_G(\phi(s))$ for $s\in H$.
As $\theta'$ is a homomorphism, we see that for $a\in A(G)$ and $s,t\in H$,
\begin{align*}
a(\phi(st)) &= \theta(a)(st) = \ip{\lambda_H(st)}{\theta(a)}
= \ip{\theta'(\lambda_H(st))}{a} = \ip{(\theta'\lambda_H(s))(\theta'\lambda_H(t))}{a} \\
&= \ip{\lambda_G(\phi(s)) \lambda_G(\phi(t))}{a}
= \ip{\lambda_G(\phi(s)\phi(t))}{a} = a( \phi(s)\phi(t) ). \end{align*}
Hence $\phi$ is a group homomorphism.

In particular, for $s\in H$,
\[ \theta'(\lambda_H(s)^*) = \theta'(\lambda_H(s^{-1})) = \lambda_G(\phi(s^{-1}))
= \lambda_G(\phi(s)^{-1}) = \lambda_G(\phi(s))^* = ( \theta'(\lambda_H(s)) )^*. \]
As $\{ \lambda_H(s) : s\in H \}$ generates $VN(H)$, we see that $\theta'$ is
a $*$-homomorphism.  Hence $\theta'$ and $\theta$ are isometries.  By
Walter's Theorem, \cite{walter}, there hence exists
\begin{enumerate}
\item either a topological group isomorphism $\psi:H\rightarrow G$, or
a topological group anti-isomorphism $\psi:H\rightarrow G$, and
\item a fixed $t_0\in G$,
\end{enumerate}
such that
\[ \theta(a)(s) = a(t_0\psi(s)) \qquad (a\in A(G), s\in H). \]
Hence $t_0\psi(s) = \phi(s)$ for all $s\in H$, and so we conclude that
$t_0 = e_G$ and $\psi = \phi$, showing that $\phi$ is bicontinuous, as required.
\end{proof}

Our main objective in this section is the following result, which shows that the canonical predual is the only isometrically induced weak$^*$-topology for $A(G)$ which respects the Hopf algebra structure. Unlike the analogous result Theorem \ref{measure.hopf}, we only consider isometric preduals.  In the previous section, 
we used Lemma~\ref{cstar_bidual_is_cstar} to allow us to reduce to the isometric case.  However, this argument seems to have no non-commutative generalisation.  The best result we can use is hence Theorem~\ref{autoc*}, which tells us that an isometric Hopf-algebra predual is automatically a C$^*$-algebra predual.

\begin{theorem}\label{fourier.hopf}
Let $G$ be a separable compact  group. If $E\subset VN(G)$ is an isometric predual for $A(G)$, making the multiplication separately weak$^*$-continuous and the coassociative multiplication continuous, then $E=C^*(G)$.
\end{theorem}

The arguments used in Section~\ref{measure_case} relied heavily
upon the Gelfand transform for a commutative C$^*$-algebra, and in
particular upon the character space.  A suitable (for us) non-commutative
analogue is the \emph{spectrum} of a C$^*$-algebra, which we now recall.
We follow the presentation in \cite[Chapter~4]{pedersen}, see also
\cite[Chapter~9]{dixmier}.  Many of these results follow fairly easily
from work of Fell in \cite{Fell} which nicely exhibits the link with
group algebras.

Let $E$ be a C$^*$-algebra, and let $P(E)$ be the set of pure states
on $E$, with the weak$^*$-topology, which need not be compact. However, when $E$ is separable, then $P(E)$ is a Polish space, that is, a separable complete metrisable space, \cite[Proposition~4.3.2]{pedersen}.
Recall that a state $\phi$ is pure if and only if the GNS representation
$\pi_\phi$ associated to $\phi$ is irreducible.  Conversely, if $\pi:
E\rightarrow\mc B(H)$ is an irreducible representation, then for any
unit vector $u\in H$, the vector state $\omega(\pi;u,u)$ is a pure state.
%Pure states $\phi,\psi\in P(E)$ give rise to equivalent representations
%$\pi_\phi$ and $\pi_\psi$ if, and only if, $\phi = u^*\cdot\psi\cdot u$
%for some unitary $u\in E$, \cite[Proposition~3.13.4]{pedersen}.

Let $\check E$ be the Primitive Ideal Space of $E$, that is, the
collection of kernels of irreducible representations.  We give $\check E$
the Hull-Kernel topology, so that the closed sets are of the form
\[ \{ t\in\check E : I\subseteq t \} \qquad (I\subseteq E). \]
Alternatively, the topology can be defined by the observation that the map
$P(E)\rightarrow\check E; \phi\mapsto\ker\pi_\phi$ is open and continuous, 
\cite[Theorem~4.3.3]{pedersen}.  Furthermore, $\check E$ is a Baire space and
$\check E$ is locally compact, but not necessarily Hausdorff,
\cite[Theorem~4.3.5, Proposition~4.4.4]{pedersen}.
%For an integer $n>0$, let ${}_n\check E \subseteq \check E$ be the primitive
%ideals corresponding to irreducible representations of dimension at most $n$.
%Let $\check E_n = {}_n\check E\setminus{}_{n-1}\check E$.  Then, ${}_n\check E$
%is closed and $\check E_n$ is Hausdorff in its relative
%topology, \cite[Proposition~4.4.10]{pedersen}.

Let $\hat E$ be the equivalence classes of irreducible representations
of $E$, so that there is a natural surjection $\hat E\rightarrow\check E$,
which we use to induce a topology on $\hat E$, by defining this map
to be open and continuous.  In many ways $\check E$ is
easier to deal with that $\hat E$, but $\hat E$ carries more information.
Fortunately, in our case, we have no problem, as we shall be dealing with
\emph{Type~I} C$^*$-algebras.  Recall that a C$^*$-algebra $E$ is
of Type~I when every irreducible representation $\pi:E\rightarrow\mc B(H)$ satisfies
$\mc K(H)\subseteq \pi(E)$. Also $E$ is a Type I if, and only
if, $\hat E = \check E$, \cite[Theorem~9.1]{dixmier}.

%Finally, we come to the key tool for our purposes.  Let $x\in E$ be positive, and define
%$\check x:\hat E\rightarrow\mathbb R$ by $\check x(\pi) = \|\pi(x)\|$ for $\pi\in\hat E$.
%Actually, $\check x(\pi) = \|x + \ker\pi \|_{E/\ker\pi}$, so that $\check x$ is a well-defined
%function on $\check E$. Then $\check x$ is continuous on $\check E$ if and only if $\check E$
%is Hausdorff, \cite[Proposition~4.4.5]{pedersen}.  Suppose now that $\check E$ is Hausdorff.
%By taking the complexification, we get a map $E\rightarrow C(\check E); x\mapsto \check x$.
%Let $\mc M(E)$ be the \emph{multiplier} algebra of $E$, so that $\mc M(E)$ can be identified
%with the C$^*$-subalgebra
%\[ \mc M(E) = \{ x\in E'' : xy,yx\in E\ ( y\in E ) \} \subseteq E''. \]
%Then the map $x\mapsto\check x$ extends to a map $\mc M(E)\rightarrow C(\check E)$ and the
%Dauns-Hofmann Theorem, \cite[Corollary~4.4.8]{pedersen}, tells us that the restriction of
%this map is an isomorphism of algebras between the centre of $\mc M(E)$ and $C(\check E)$.

Let $G$ be a locally compact group.  As the unitary representations of $G$
and $C^*(G)$ agree, for our purposes, we may define the Fell topology on
$\hat G$ to be the topology on $C^*(G)\hat{}$.  It is well known that when $G$
is a compact group, then $\hat G$ has the discrete topology and that each
member of $\hat G$ is finite dimensional, \cite[Theorem~15.1.3]{dixmier}.
For $\pi\in\hat G$, let $\pi:G\rightarrow\mc B(H_\pi)$ where $H_\pi$ is a
$\dim(\pi)$-dimensional Hilbert space.  It hence follows that
%\[ C^*(G) \cong c_0-\bigoplus_{\pi\in\hat G} \mathbb M_{\dim(\pi)},
%\quad VN(G) \cong \ell^\infty-\bigoplus_{\pi\in\hat G} \mathbb M_{\dim(\pi)}. \]
\[ C^*(G) \cong c_0-\bigoplus_{\pi\in\hat G} \mc B(H_\pi),
\quad VN(G) \cong \ell^\infty-\bigoplus_{\pi\in\hat G} \mc B(H_\pi). \]
Let $\mc T(H_\pi)$ be the trace-class operators on $\mc B(H_\pi)$, so as $H_\pi$ is
finite-dimensional, we have that $\mc B(H_\pi)' = \mc T(H_\pi)$ and
$\mc T(H_\pi)' = \mc B(H_\pi)$.
%Here $\mathbb M_n$ is the space of complex $n\times n$ matrices.  Let $\mathbb T_n$
%be its dual space, the space of complex $n\times n$ matrices with the trace-class norm.
Then, as a Banach space,
%\[ A(G) \cong \ell^1-\bigoplus_{\pi\in\hat G} \mathbb T_{\dim(\pi)}, \]
\[ A(G) \cong \ell^1-\bigoplus_{\pi\in\hat G} \mc T(H_\pi), \]
but the algebra product on $A(G)$ is not easily expressed under this identification.
Essentially, this is because the tensor product of two irreducible representations need not be irreducible.

For the remainder of this section, we shall let $E\subseteq VN(G)$ be an isometric Hopf
algebra predual for $A(G)$.  Hence by Proposition~\ref{fourier_hopf=>alg}, $E$ is a
subalgebra of $VN(G)$ and so by Theorem \ref{autoc*} $E$ is a C$^*$-algebra.
By Propositions~\ref{states} and~\ref{reps}, every irreducible representation
of $E$ is induced in the canonical way by an irreducible representation
of $C^*(G)$.  That is, as sets, $\hat E = \hat G$, although the
topologies may differ.  Thus $E$ is a Type~I C$^*$-algebra, and so
$\check E = \hat E$.  The above direct sum matrix form for $C^*(G)$ 
and $VN(G)$ shows that the irreducible representations of $\hat G$
are simply the projections onto one of the factors $\mc B(H_\pi)$.  %$\mathbb M_{\dim(\pi)}$.
Just as in section~\ref{measure_case}, it will suffice to show that $\hat{E}$ has the same topology as $\hat{G}$.

Let us note that singletons in $\hat E$ are closed.
Indeed, let $\pi\in\hat E$, and suppose that $\phi\in\hat E$ is in the closure of
$\{\pi\}$.  As $\hat E = \check E$, using the hull-kernel topology, we see that
$\ker\pi \subseteq \ker\phi$.  Hence we get a natural map $E / \ker\pi \rightarrow
E / \ker\phi$.  As $\pi:E\rightarrow\mc B(H_\pi)$ is irreducible, with $H_\pi$ being
finite-dimensional, $E/\ker\pi$ is simple, and so either this
map is the zero map, so $\ker\phi=E$, that is, $\phi=0$, contradiction; or the map
is injective, and hence $\ker\pi=\ker\phi$, so $\pi$ and $\phi$ are equivalent,
that is, $\pi=\phi$ in $\hat E$.

\begin{proposition}\label{discrete_enough}
With notation as above, $E=C^*(G)$ if, and only if, $\hat E$ has
the discrete topology.
\end{proposition}
\begin{proof}
Clearly we need only show that if $\hat E$ is discrete, then $E=C^*(G)$.
%This follows almost immediately from the Dauns-Hofmann Theorem.
%As $\hat E$ is discrete, it is Hausdorff, and for any $\pi\in\hat E$,
%the indicator function of $\{\pi\}$ is a continuous function.  Hence
%this function corresponds to an element $p$ in the centre of $\mc M(E)$.
%It is clear that $p$ is the projection corresponding to taking
%$x\in E$ to the component of $x$ in the matrix algebra
%corresponding to $\pi$.  As $p$ is a multiplier and $\pi$ was arbitrary,
%we immediately see that $E = C^*(G)$, as required.
As $\hat E$ is discrete, for any $\pi\in\hat E$, the singleton $\{\pi\}$ is open (and
closed by a comment above).  %Define $z=(z_\rho)_{\rho\in\hat G}\in VN(G)$ by
%$z_\rho = 1$ if $\rho=\pi$, and $z_\rho=0$ otherwise.  We claim that $z$ is a
%\emph{multiplier} of $E$, that is, $zx\in E$ for all $x=(x_\rho) \in E$.

By the hull-kernel topology, as $\check E\setminus \{\pi\}$ is closed, we have
that $\ker(\check E\setminus\{\pi\}) \not\subseteq \ker\pi$, where
\[ \ker( \check E \setminus \{\pi\} ) = \big\{ w=(w_\rho)\in E :
w_\rho=0 \ (\rho\not=\pi) \big\}. \]
So there exists $w=(w_\rho)\in E$ with $w_\rho=0$ for all $\rho\not=\pi$, and
with $w_\pi\not=0$.  As $\pi$ is an irreducible representation of $E$, for any
$a\in\mc B(H_\pi)$, we can find $y\in E$ with $y_\pi=a$.  As $\mc B(H_\pi)$ is
simple, it follows that for any $a\in\mc B(H_\pi)$, we can find $\sum_n y^{(1)}_n
\otimes y^{(2)}_n \in E\otimes E$ such that, if $z = \sum_n y^{(1)}_n w y^{(2)}_n$,
then $z_\pi=a$; clearly we have that $z_\rho=0$ for $\rho\not=\pi$.

It is now immediate that $E$ has the same form as $C^*(G)$, so $E = C^*(G)$
as subspaces of $VN(G)$.
\end{proof}

We now sketch some theory about compact groups and their representations.
See, for example, \cite[Chapter~15]{dixmier}.  Let $G$ be a compact group, and let $\pi:G\rightarrow\mc B(H)$ be a
finite-dimensional representation.  The \emph{character}
of $\pi$ is the map $\chi_\pi:G\rightarrow\mathbb C$ defined by
\[ \chi_\pi(g) = \operatorname{Tr}(\pi(g)) \qquad (g\in G). \]
A \emph{class function} on $G$ is some function constant on conjugate
classes.  Then $\chi_\pi$ is a continuous class function, and $\chi_\pi$
only depends upon the equivalence class of $\pi$.  The collection
$\{ \chi_\pi : \pi\in\hat G \}$ is dense in the space of continuous class
functions equipped with the supremum norm, and also forms an orthonormal basis for the space of $L^2$ class functions with the inner product
\[ [\chi_1, \chi_2] = \int_G \chi_1(g) \overline{ \chi_2(g) }
\ dg. \]

The character of a finite-dimensional representation $\pi$ determines
the equivalence class of $\pi$ in the following way.  For $n\in\mathbb N$
and $\rho\in\hat G$, write $n\rho$ for the representation
$\rho\oplus\cdots\oplus\rho$, where $\rho$ is repeated $n$ times.
Then $\pi$ is equivalent to
\[ \sum_{\rho\in\hat G} \ [\chi_\pi, \chi_\rho] \rho. \]
Here we really only sum over $\rho\in\hat G$ such that
$[\chi_\pi, \chi_\rho]\not=0$, and it is part of the theory that
$[\chi_\pi, \chi_\rho]$ is always a positive integer. A simple calculation shows that $\chi_{\pi\otimes\rho} 
= \chi_\pi \chi_\rho$, so characters allow us to work out the
equivalence class of $\pi\otimes\rho$.  However, we have no concrete
way to actually find a unitary which implements this equivalence. Recall too, that the contragradient representation $\pi^*$ associated to a representation $\pi$ of $G$ on $H$ is given by representing $G$ on the conjugate Hilbert space $\overline{H}$. The character of $\pi^*$ is given by
\[ \chi_{\pi^*}(g) = \overline{\chi_{\pi}(g)} \qquad (g\in G). \]

For $\pi_0\in\hat G$, let $\alpha,\beta\in H_{\pi_0}$, and let $\alpha\otimes
\beta\in A(G)$ be the (normal) functional on $VN(G)$ defined by
\[ \ip{x}{\alpha\otimes\beta} = (x_{\pi_0}(\alpha)|\beta)
\qquad \big( x=(x_\pi) \in VN(G) = \ell^\infty( \mc B(H_\pi) ) \big). \]

\begin{lemma}\label{fourier.tech1}
Let $G$ be a compact group, let $\pi_0,\rho,\pi\in\hat G$ be such
that $[\chi_{\pi_0} \chi_\rho, \chi_\pi]>0$.  Then there exists $\xi\in
H_{\pi_0}$ with $\|\xi\|=1$, $\alpha\in H_\rho$ and $w\in\mc B(H_\pi)$
such that, if $a=\xi\otimes\xi, b=\alpha\otimes\alpha\in A(G)$, and
$c = ab = (c_\pi) \in \ell^1(\mc T(H_\pi))$, then $[w,c_\pi]\not=0$.
\end{lemma}
\begin{proof}
We know that the representation $\pi_0 \otimes \rho$ is equivalent
to the representation
\[ \sum_{\phi\in\hat G} \ [\chi_{\pi_0} \chi_\rho , \chi_\phi] \phi. \]
Let $I = \{ \phi\in\hat G : [\chi_{\pi_0} \chi_\rho , \chi_\phi]>0 \}$
a finite set and let $J$ be the finite collection of irreducible
representations of $G$ formed by taking $\phi\in I$ with the multiplicity
$[\chi_{\pi_0} \chi_\rho , \chi_\phi]$.  Hence there exists some
unitary
\[ U:H_{\pi_0} \otimes H_\rho \rightarrow \bigoplus_{\phi\in J} H_\phi, \]
such that, if $U_\phi$ is the component of $U$ mapping to $H_\phi$ for
each $\phi\in J$, then
\[ U_\phi (\pi_0(g) \otimes \rho(g)) = \phi(g) U_\phi
\qquad (\phi\in J, g\in G). \]
Notice that each $U_\phi$ is a partial isometry.  With $a$ and $b$ as
defined, we have
\[ \ip{x}{c} = \ip{x}{ab} = \sum_{\phi\in J} (x_\phi U_\phi(\xi\otimes\alpha)|
U_\phi(\xi\otimes\alpha)) \qquad (x\in VN(G)). \]
So if our claim is false, as we can vary $w$, we must have that
$U_\pi(\xi\otimes\alpha)=0$, for all $\xi$ and $\alpha$, which implies
that $U_\pi=0$.  As $[\chi_{\pi_0} \chi_\rho, \chi_\pi]>0$, we have that
$\pi\in J$, giving a contradiction, as required.
\end{proof}

We are now in a position to establish Theorem \ref{fourier.hopf}.
\begin{proof}[Proof of Theorem \ref{fourier.hopf}.]
By Proposition~\ref{discrete_enough}, we need to show that $\hat E$
is discrete.  
Suppose, towards a contradiction, that $\hat E$ is not discrete,
so we can find $\pi_0\in\hat E$ with $\{\pi_0\}$ not open.
As each singleton in $\hat E$ is closed, and $\hat E$ is a countable Baire space, there exists $\pi_1\in\hat E$ with $\{\pi_1\}$ open.
We claim that there exists some $\rho\in\hat E$ with
$[\chi_{\pi_0} \chi_\rho , \chi_{\pi_1}]\not=0$.  Notice that
\[ [\chi_{\pi_0} \chi_\rho , \chi_{\pi_1}]
= [\chi_\rho , \chi_{\pi_0^*} \chi_{\pi_1}], \]
so if our claim is false, then $\chi_\rho$ is orthogonal to
$\chi_{\pi_0^*}\chi_{\pi_1}$ for every $\rho\in\hat E$.  However,
$\{ \chi_\rho : \rho\in\hat E \}$ is an orthonormal basis of the $L^2$
class functions, so $\chi_{\pi_0^*}\chi_{\pi_1} = 0$
in $L^2(G)$.  As $\chi_{\pi_0^*}\chi_{\pi_1}$ is continuous, this
implies that $\chi_{\pi_0^*}\chi_{\pi_1}=0$ identically.  This is
patently untrue, simply evaluate at the identity of $G$. Thus we can find some $\rho$, as claimed.

The set $I = \{ \pi\in\hat E : [\chi_\pi \chi_\rho , \chi_{\pi_1}]\not=0 \}$
is finite and contains $\pi_0$.  As singletons in $\hat E$ are closed, it follows
that finite sets are closed, and so we can find an open set $U\subseteq\hat E$ such
that $U\cap I = \{\pi_0\}$.
%By Lemma \ref{fourier.tech2}, there exists an open
%set $U\subseteq\hat E$ which contains $\pi_0$, and such that $U\cap I = \{\pi_0\}$.
Let $f:P(E)\rightarrow\hat E$ be the natural map, which is open and
continuous.  We see that for $\pi\in\hat G$,
\[ f^{-1}(\{\pi\}) = \big\{ \xi\otimes\xi : \xi\in H_\pi, \|\xi\|=1 \}. \]
That $\{\pi_0\}$ is not open means that $f^{-1}(\{\pi_0\})$ is not open
in $P(E)$, that is, there exists some $\xi_0\in H_{\pi_0}$ with $\|\xi_0\|=1$,
such that for each finite set $F\subseteq E$ and $\epsilon>0$,
there exists $\pi\not=\pi_0$ and $\eta\in H_\pi$ with $\|\eta\|=1$ and
\[ \big| \ip{\xi_0\otimes\xi_0 - \eta\otimes\eta}{x} \big| < \epsilon
\qquad (x\in F). \]
That is, every weak$^*$-open neighbourhood of $\xi_0\otimes\xi_0$ contains
some member of $f^{-1}(\{\pi\})$ for some $\pi\not=\pi_0$.

As $f^{-1}(U)$ is open in $P(E)$ and contains $\xi_0\otimes\xi_0$, there
exists a finite set $F_0\subseteq E$ such that, if $\pi\in\hat G$,
$\eta\in H_\pi$ with $\|\eta\|=1$ and
\[ \big| \ip{\xi_0\otimes\xi_0 - \eta\otimes\eta}{x} \big| < 1
\qquad (x\in F_0), \]
then $\pi\in U$.  Let $F\subseteq E$ be a finite set, so by the
previous paragraph, there exists $\pi_F\not=\pi_0$ and
$\eta_F\in H_{\pi_F}$ with $\|\eta_F\|=1$ and
\[ \big| \ip{\xi_0\otimes\xi_0 - \eta_F\otimes\eta_F}{x} \big| < |F|^{-1}
\qquad (x\in F\cup F_0). \]
In particular, $\pi_F\in U$.  Let $a_0 = \xi_0\otimes\xi_0\in A(G)$ and
for each $F\subseteq E$ finite, let $a_F=\eta_F\otimes\eta_F\in A(G)$.
By construction, $a_F \rightarrow a_0$ in the weak$^*$-topology
on $A(G)$ induced by $E$.

Now suppose that $\xi\in H_{\pi_0}$ is arbitrary, with $\|\xi\|=1$.  As the
GNS representations for $\xi_0$ and $\xi$ are equivalent,
there exists some unitary $u=(u_{\pi})\in E$ with $u_{\pi_0}(\xi_0)=\xi$.
Then, for $x\in E$,
\[ \lim_F \ip{ua_Fu^*}{x} = \lim_F \ip{a_F}{u^*xu}
= \ip{a_0}{u^*xu} = (u_{\pi_0}^* x_{\pi_0} u_{\pi_0}(\xi_0)|\xi_0)
= (x_{\pi_0}(\xi)|\xi), \]
so that $ua_Fu^* \rightarrow \xi\otimes\xi$ weak$^*$.

%As $\{\pi_1\}$ is open in $\hat E$, by the Dauns-Hoffmann theorem,
%for any $w\in\mc B(H_{\pi_1})$, there exists $x=(x_\pi)\in E$ with
%$x_\pi = 0$ for $\pi\not=\pi_1$, and with $x_{\pi_1} = w$.
As $\{\pi_1\}$ is open in $\hat E$, by the same argument as used in the proof
of Proposition~\ref{discrete_enough} we see that for any $w\in\mc B(H_{\pi_1})$, there exists
$x=(x_\pi)\in E$ with $x_\pi = 0$ for $\pi\not=\pi_1$, and with $x_{\pi_1} = w$.
Let $\rho\in\hat G$.
Let $\alpha\in H_\rho$, let $b = \alpha\otimes\alpha\in A(G)$,
let $\xi\in H_{\pi_0}$ with $\|\xi\|=1$, and let $a=\xi\otimes\xi\in
A(G)$.  Then let $c=ab\in A(G)$, say $c = (c_\pi) \in \ell^1-\bigoplus
\mathcal T(H_\pi)$.  We see that
\[ \ip{c}{x} = \ip{c_{\pi_1}}{w} \not=0 \]
for some choice of $\alpha, \xi$ and $w$, by Lemma \ref{fourier.tech1}.

Let $u\in E$ be some unitary with $u_{\pi_0}(\xi_0) = \xi$.
Let $c^{(F)} = (ua_Fu^*)b$ for $F\subseteq E$ finite, and
suppose that $c^{(F)}_{\pi_1} \not= 0$.  Thus $[\chi_\rho \chi_{\pi_F},
\chi_{\pi_1}]>0$, but as $\pi_F\in U$, this is a contradiction.
In conclusion, we have that $\ip{ua_Fu^*}{b\cdot x}=0$ for
each finite $F\subseteq E$, so that
\[ \ip{c}{x} = \ip{a}{b\cdot x} = \lim_F \ip{u a_F u^*}{b\cdot x} = 0, \]
a contradiction.  Hence $\hat E$ is discrete, and the result follows from Proposition \ref{discrete_enough}.
\end{proof}

\section{Algebras with unique preduals}\label{Uniqueness}
For certain classes of dual Banach algebras, the predual is uniquely determined so that there is one weak$^*$-topology.  The first example of this phenomena is in \cite[Theorem 4.4]{Daws}, where $\mathcal B(E)$ is shown to have a unique predual for any reflexive Banach space $E$ with the approximation property.  A careful examination of the proof of this theorem yields the following result, since the hypothesis stated below are the only properties of $\mathcal B(E)$ used.
We refer the reader to the discussion after \cite[Theorem~2.6.15]{Dales} for details on the Arens products.

\begin{theorem}\label{Uniqueness.Arens}
Let $\mc A$ be an Arens regular Banach algebra such that $\mc A''$ is
unital, and $\mc A$ is an ideal in $\mc A''$.  Then $\mc A'$ is the unique
predual of $\mc A''$.
\end{theorem}

Von Neumann algebras can be characterised abstractly as those $C^*$-algebras which are isometric to the dual space of some Banach space.  There is a unique \emph{isometric} weak$^*$-topology on a von Neumann algebra and this topology make the multiplication separately continuous and the adjoint continuous. In contrast, a classical example of Pelcyznski \cite{Pel} shows that the commutative non-isomorphic von Neumann algebras $\ell^\infty$ and $L^\infty[0,1]$ are isomorphic as Banach spaces, and so the non-isomorphic Banach spaces $\ell^1$ and $L^1[0,1]$ induce two distinct weak$^*$-topologies on $\ell^\infty$ --- of course the topology induced by $L^1[0,1]$ does not respect the von Neumann algebra structure. It was shown in
\cite{Daws} that if $\theta:\mc M\rightarrow\mc N$ is merely a Banach
algebra isomorphism, and $\mc M$ and $\mc N$ are commutative von
Neumann algebras, then $\theta$ is weak$^*$-continuous and so a commutative von Neumann algebra has a unique weak$^*$-topology making the multiplication separately continuous.  The weak$^*$-continuity of the adjoint follows for free.  The theorem below extends this to the non-commutative setting by passing through maximal abelian subalgebras.

\begin{theorem}\label{Uniqueness.VNA}
Let $\mc M$ be a von Neumann algebra, let $\mc B$ be a dual Banach
algebra with predual $\mc B_*$, and let $\theta:\mc M\rightarrow\mc B$
be a Banach algebra isomorphism.  Then $\theta$ is weak$^*$-continuous.
\end{theorem}
\begin{proof}
As in the proof of \cite[Theorem~5.1]{Daws}, it is enough to show
that if $X\subseteq\mc M'$ is a predual, then $X$ is the usual predual
for a von Neumann algebra, that is, each functional in $X$ is normal.
This equivalence follows by setting $X=\theta'(\mc B_*)$.

So let $X\subseteq\mc M'$ be a predual, and pick $\mu\in X$.  Then,
by \cite[Chapter~III, Corollary~3.11]{tak}, $\mu$ is normal if, and only
if, the restriction of $\mu$ to each maximal abelian self-adjoint subalgebra (masa) is normal.
Take a masa $\mc A$ in $\mc M$. By maximality, $\mc A$ is weak$^*$-closed in any weak$^*$-topology
arising a predual, and in particular in that induced by $X$. Identify $\mc M$ with $X'$,
so that $\mc A$ has the predual $X / {^\perp}\mc A$.  By \cite[Theorem~5.1]{Daws}, it follows
that each member of $X / {^\perp}\mc A$ is normal.  Hence $\mu+{^\perp}\mc A$ is normal,
and $\ip{a}{\mu+{^\perp}\mc A} = \ip{a}{\mu}$ for $a\in\mc A$, so that
$\mu$, when restricted to $\mc A$, is normal.  We hence conclude that
$\mu$ is normal, as required.
\end{proof}

\begin{corollary}
There is a unique weak$^*$-topology on a von Neumann algebra which makes the multiplication map separately continuous.
\end{corollary}

\bigskip

\begin{tabular*}{\textwidth}{l@{\hspace*{.5cm}}l@{\hspace*{.5cm}}l}
Matthew Daws&Hung Le Pham&Stuart White\\
School of Mathematics &Department of Mathematical&Department of Mathematics\\
&and Statistical Sciences&\\
University of Leeds&University of Alberta           &University of Glasgow\\
Leeds        &Edmonton  &Glasgow\\
LS2 9JT&T6G 2E1&G12 8QW\\
UK&Canada &UK\\
[1ex]
\texttt{matt.daws@cantab.net}&\texttt{hlpham@math.ualberta.ca}&\texttt{s.white@maths.gla.ac.uk}
\end{tabular*}


\begin{thebibliography}{9}

\newcommand{\bibbook}[3]{\textsc{#1}, \emph{#2}, (#3).}
\newcommand{\bibpaper}[6]{\textsc{#1}, `#2', \emph{#3} #4 (#5) #6.}
\newcommand{\bibpreprint}[2]{\textsc{#1}, `#2', preprint.}

\bibitem{TB} \bibpaper{L. Baggett, K.\,F. Taylor}
   {Groups with completely reducible regular representation}
   {Proc. Amer. Math. Soc.}{72}{1978}{593--600}

\bibitem{Dales} \bibbook{H.\,G. Dales}
   {Banach algebras and automatic continuity}
   {Clarendon Press, Oxford, 2000}

\bibitem{Daws} \bibpaper{M. Daws}
   {Dual Banach algebras: representations and injectivity}
   {Studia Math.}{178}{2007}{231--275}
  
  \bibitem{DHSW} \textsc{M. Daws, R. Haydon, T. Schulmprecht, S. White}, `Shift invariant preduals', Manuscript in preparation.

\bibitem{dixmier} \bibbook{J. Dixmier}
   {$C\sp*$-algebras. Translated from the French by Francis Jellett}
   {North-Holland Mathematical Library, Vol. 15. North-Holland Publishing Co., Amsterdam-New York-Oxford, 1977}

\bibitem{ER} \bibpaper{E.\,G. Effros, Z.-J. Ruan}
   {Operator space tensor products and Hopf convolution algebras}
   {J. Operator Theory}{50}{2003}{131--156}

\bibitem{ellis} \bibpaper{R. Ellis}
   {Locally compact transformation groups}
   {Duke Math. J.}{24}{1957}{119--125}

\bibitem{Eymard} \bibpaper{P. Eymard}{L'algebr\'a de Fourier d'un groupe
   localement compact}{Bull. Soc. Math. France}{92}{1962}{181--236}

\bibitem{Fell} \bibpaper{J.\,M.\,G. Fell}
   {The dual spaces of $C\sp{*} $-algebras}
   {Trans. Amer. Math. Soc.}{94}{1960}{365--403}

\bibitem{johnson} \bibpaper{B.\,E. Johnson}
   {Isometric isomorphisms of measure algebras}
   {Proc. Amer. Math. Soc.}{15}{1964}{186--188}

\bibitem{palmer} \bibpaper{T.\,W. Palmer}
   {Arens multiplication and a characterization of $w\sp{*} $-algebras}
   {Proc. Amer. Math. Soc.}{44}{1974}{81--87}

\bibitem{Pel} \bibpaper{Pe\l czy\'nski, A.}
   {On the isomorphism of the spaces $m$ and $M$.}
   {Bull. Acad. Polon. Sci. Sér. Sci. Math. Astr. Phys.}{6}{1958}{695--696}

\bibitem{pedersen} \bibbook{G.\,K. Pedersen}
   {$C\sp{*} $-algebras and their automorphism groups}
   {London Mathematical Society Monographs, 14. Academic Press, Inc., London-New York, 1979}

\bibitem{Quigg} \bibpaper{J.\,C. Quigg}
   {Approximately periodic functionals on $C\sp *$-algebras and von Neumann algebras}
   {Canad. J. Math.}{37}{1985}{769--784}

\bibitem{runde} \bibpaper{V. Runde}
   {Amenability for dual Banach algebras}
   {Studia Math.}{148}{2001}{47--66}

\bibitem{Ryan} \bibbook{R. Ryan}
   {Introduction to Tensor Products of Banach Spaces}
   {Springer-Verlag, London, 2002}

\bibitem{tak} \bibbook{M. Takesaki}
   {Theory of Operator Algebras I}
   {Springer-Verlag, New York, 1979}

\bibitem{taylor} \bibpaper{K.\,F.~Taylor}
   {Geometry of the Fourier algebras and locally compact groups with atomic unitary representations}
   {Math. Ann.}{262}{1983}{183--190}

\bibitem{walter} \bibpaper{M.\,E. Walter}
   {Group duality and isomorphisms of Fourier and Fourier-Stieltjes algebras from a $W\sp{*} $-algebra point of view.}
   {Bull. Amer. Math. Soc.}{76}{1970}{1321--1325}

\end{thebibliography}
\end{document}